\documentclass{amsart}

\newtheorem{thm}{Theorem}[section]
\newtheorem{lem}[thm]{Lemma}
\newtheorem{prp}[thm]{Proposition}

\theoremstyle{definition}

\theoremstyle{remark}

\numberwithin{equation}{section}

 \newcommand{\eps}{\varepsilon}
 \newcommand{\vpi}{\varphi}
 
 \newcommand{\lo}{\mathcal{L}}
 \newcommand{\D}{\mathcal{D}}
 \newcommand{\Real}{\mathbb{R}}

 \newcommand{\norm}[1]{\left\Vert#1\right\Vert}
 \newcommand{\abs}[1]{\left\vert#1\right\vert}
 \newcommand{\set}[1]{\left\{#1\right\}}
 \newcommand{\bigset}[1]{\big\{#1\big\}}
 \newcommand{\inner}[1]{\left(#1\right)}

 \newcommand{\mc}[1]{{\mathcal #1}}

 \newcommand{\bb}[1]{{\mathbb #1}}
 \newcommand{\reff}[1]{(\ref{#1})}
 
\begin{document}
\title[The Gevrey Hypoellipticity for Fokker-Planck equations]
{The Gevrey Hypoellipticity for linear and non-linear Fokker-Planck equations}
\thanks{Research  partially supported by NSFC}
\date{}



\maketitle

\begin{center}

\author{{Hua Chen}$^{1}$ \& {Wei-Xi Li}$^1$ \& {Chao-Jiang Xu}$^{1,2}$
\thanks{}\\\vspace{0.5cm}
{}$^1$School of Mathematics and Statistics, Wuhan University,\\
Wuhan 430072, China\\
email: chenhua@whu.edu.cn\\
\hspace{0.9cm}wei-xi.li@whu.edu.cn\\[2mm]
{}$^2$Universit\'{e} de Rouen, UMR 6085-CNRS, Math\'{e}matiques\\
Avenue de l'Universit\'{e}, BR.12,76801 Saint Etienne du Rouvray,
France\\
email: Chao-Jiang.Xu@univ-rouen.fr}

\end{center}

\begin{abstract}
In this paper, we study the Gevrey regularity of weak solution for a
class of linear and quasilinear Fokker-Planck equations.
\end{abstract}

\bigbreak
\section{Introduction}\label{sect1}

Recently, a lot of progress has been made on the study for the
spatially homogeneous Boltzmann equation without angular cutoff, cf.
\cite{al-1,al-2,desv-wen1,villani} and references therein, which
shows that the singularity of collision cross-section yields some
gain of regularity in the Sobolev space frame on weak solutions for
Cauchy problem. That means, this gives the $C^\infty$ regularity of
weak solution for the spatially homogeneous Boltzmann operator
without angular cutoff. The local solutions having  the Gevrey
regularity have been constructed in \cite{ukai} for initial data
having the same Gevrey regularity, and a genearal Gevrey regularity
results have given in \cite{MUXY1} for spatially homogeneous and
linear Boltzmann equation of Cauchy problem for any initial data. In
the other word, there is the smoothness effet similary to heat
equation.

However, there is no general theory for the spatially inhomogeneous
problems. It is now a kinetic equation in which the diffusion part
is nonlinear operator of velocity variable. In \cite{aumxy}, by
using the uncertainty principle and microlocal analysis, they obtain
a $C^\infty$ regularity results for linear spatially inhomogeneous
Boltzmann equation without angular cutoff.

In this paper, we will study the Gevrey regularity of weak solution
for the the following  Fokker-Planck operator in $\bb{R}^{2n+1}$
\begin{eqnarray}\label{Fokker-Planck1}
\mc{L}=\partial_t+v\cdot\partial_x-a(t,x,v)\triangle_v,
\end{eqnarray}
where $\triangle_v$ is Laplace operator of velocity variables $v$.

The motivation of study for this class of operator is , as in
\cite{MX}, attempt to study inhomogenous Boltzmann equations without
angular cutoff and non linear Vlasov-Fokker-Planck equation (see
\cite{helffer-nier,herau-nier}).

Before stating the result, we recall the definition of Gevrey class
function. Let $U$ be an open subset of $\mathbb{R}^N$ and $f$ be a
real function defined in $U$. We say $f\in G^s(U)(s\geq1)$ if $f\in
C^\infty(U)$ and for any compact subset $K$ of $U$, there exists a
constant $C=C_K$, depending only on $K$, such that for all
multi-indices $\alpha\in\mathbb{N}^N$ and for all $x\in K$
\begin{eqnarray}\label{gevrey}
|\partial^\alpha{f}(x)| \leq C_K^{|\alpha|+1}(|\alpha|!)^s.
\end{eqnarray}
Denote by $\bar{U}$ the closure of $U$ in $\bb{R}^N.$ we say $f\in
G^s(\bar{U})$ if $f\in G^s(W)$ for some open neighborhood $W$ of
$\bar{U}.$ The estimate {\rm (\ref{gevrey})} for $x\in K$ is valid
if and only if the following one is valid {\rm( cf.Chen
hua-Rodino[\ref{CR}] or Rodino[\ref{Ro}])}:
$$
\|\partial^\alpha{f}\|_{L^2(K)}\leq
C_K^{|\alpha|+1}(|\alpha|)^{s|\alpha|}.
$$
In this paper, we use the above estimate in $L^2$.

We say an operator $P$ is $G^s$ hypoelliptic in $U$ if $u\in \D',
Pu\in G^s(U)$ implies $u\in G^s(U).$ Likewise, we say an operator
$P$ is $C^\infty$ hypoelliptic in $U$ if $u\in \D', Pu\in
C^\infty(U) $ implies $u\in C^\infty(U).$

The operator $\mc{L}$ satisfies the H\"{o}rmander' condition.~~By
virtue of the results of H\"{o}rmander [5], we know that $\mc{L}$ is
$C^\infty$ hypoelliptic. In the aspect  of Gevrey class,
Derridj-Zuily [\ref{DZ}] proved that $\mc{L}$ is $G^s$-hypoelliptic
for $s>6$ in a general form of H\"ormander's operators.

In this paper, we improve firstly the results of [\ref{DZ}] for
Fokker-Planck operators as the following theorems.

\begin{thm}\label{th1}
For any $s\geq 3$,  if the coefficient $ a $ is in
$G^s(\mathbb{R}^{2n+1})$  and
 $a(t,x,v)\geq c_0>0,$  then the operator $\mc{L}$ given in
(\ref{Fokker-Planck1}) is $G^s$ hypoelliptic in $\mathbb{R}^{2n+1}$,

\end{thm}

Of course, Theorem \ref{th1} is also true for the following general
operators,
\[
\tilde{\mc{L}}=\partial_t+A(v)\cdot\partial_x-\sum^n_{j, k=1} a_{j
k}(t,x, v)
\partial^2_{v_j v_k}
\]
in an open domain $U$ of ${\bb R^{2n+1}}$, where $A$ is a non
singular  $n\times n$ constant matrix, $\big(a_{j k}(t,x, v)\big)$
is positive defined on $U$ and belongs to $G^s(U).$\\

\noindent {\bf Remark}  Our results is a local and interior
regularity results, that means if there exists a weak solution in
$\D'$, then this solution is in Gevrey class in interior of domain.
So that if the weak solution is a solution of the Cauchy problem, we
don't need the regularity of initial data.

Secondly, we consider the quasi-linear equation
\begin{equation}\label{++1.2}
\partial_t u+v \cdot \nabla _x u-a {\triangle_v}u=F(t, x, v, u,
\nabla_v u)
\end{equation}
where $F$ is nonlinear function of real variable $(t, x, v, w, p)$.

\begin{thm}\label{th3}
Let $u$ be a weak solution of equation (\ref{++1.2}) such that $u,
\nabla_v u \in L^{\infty}_{loc}(\bb R^{2n+1})$, then
\[
u\in  G^s(\mathbb{R}^{2n+1})
\]
for any $s\geq3$, provided the coefficients $ a $ is in
$G^s(\mathbb{R}^{2n+1})$,  $a(t,x,v)\geq c_0>0$ and nonlinear
function $ F(t, x, v, w,p) $ is in $G^s(\mathbb{R}^{2n+2+n})$.
\end{thm}
\noindent{\bf Remark} : If the nonlinear term $F$ is independent of
$p$ or in the form of
$$
\nabla_v (F(t, x, v, u)),
$$
we can suppose that the weak solution $u \in L^{\infty}_{loc}(\bb
R^{2n+1})$.

The plan of this paper is as follows : In  section \ref{sect2}, we
obtain a sharp subelliptic estimate for the Fokker-Planck operator
$\mc L$ via direct computation, and then prove the Gevrey
hypoellipticity of $\mc L$. In section \ref{sect3}, we prove the
Gevrey regularity for the weak solutions of the quasi-linear
Fokker-Planck equation.

\section{Subelliptic estimate}
\label{sect2} \setcounter{equation}{0}

We recall firstly some notations, $\|\cdot\|_\kappa,
\kappa\in\mathbb R,$
 is  the classical Sobolev norm in $H^\kappa(\mathbb R^{2n+1})$, and
 $(h,~k)$ is the inner product of $h, k\in L^2(\bb R^{2n+1})$. Moreover
 if $f,g\in C_0^\infty(\bb R^{2n+1}),$ from H\"{o}lder inequality and Young
inequality, for any $\varepsilon>0,$
\begin{eqnarray}\label{Young}
|(f,~g)|\leq \|h\|_\kappa \|g\|_{-\kappa}\leq
\frac{\varepsilon\|h\|_\kappa^2}{2}+\frac{\|g\|_{-\kappa}^2}{2\varepsilon}.
\end{eqnarray}

We have also the interpolation inequality for Sobolev space, for any
$\varepsilon>0$ and any $r_1<r_2<r_3,$
\begin{eqnarray}\label{interpolation}
\|h\|_{r_2}\leq \varepsilon
\|h\|_{r_3}+\varepsilon^{-(r_2-r_1)/(r_3-r_2)}\|h\|_{r_1}.
\end{eqnarray}

Let $\Omega$ be an open subset of $\bb R^{2n+1}$.  We denote by
$S^m=S^m(\Omega), m\in \bb R,$ the symbol space of classical
pseudo-differential operator and $P=P(t,x,v,D_t,D_x,D_v)\in{\rm
Op}(S^m)$ a pseudo-differential operator of symbol
$p(t,x,v;\tau,\xi,\eta)\in S^m.$ If $P\in{\rm Op}(S^m)$, then $P$ is
a continuous operator from $H_{c}^\kappa(\Omega)$ to
$H_{loc}^{\kappa-m}(\Omega)$. Here $H_{c}^\kappa(\Omega)$ is the
subspace of $H^\kappa(R^{2n+1})$ consisting of  the distributions
having their compact support in $\Omega$, and
$H_{loc}^{\kappa-m}(\Omega)$ consists of the distributions $h$ such
that $\phi h\in H^{\kappa-m}(\bb R^{2n+1})$ for any $\phi\in
C_0^\infty(\Omega)$. The more properties can be found in the Treves'
book [\ref{Treves}]. Remark that if $P_1\in {\rm Op}(S^{m_1})$,
$P_2\in {\rm Op}(S^{m_2})$, then $[P_1,~P_2]\in {\rm
Op}(S^{m_1+m_2-1}).$

\vspace{0.3cm} Now we show a sharp subelliptic estimate for the
operator $\mc L$, our proof bases on the work of Bouchut \cite{Bou}
and Morimoto-Xu \cite{MX}.

\begin{prp}\label{prp2.1}
Let $K$ be a compact subset of $\bb{R}^{2n+1}.$  Then for any $r\geq
0,$ there exists a constant $C_{K,r},$ depending only on $K$ and
$r$, such that for any $f\in C_0^\infty(K),$
\begin{eqnarray}\label{subestimate}
\|f\|_r\leq C_{K,r}\{~\|\mc{L}f\|_{r-2/3}+\|f\|_0~\}.
\end{eqnarray}

\end{prp}

To simplify the notation, in this section we will denote by $C_K$
the different suitable constants depending only on $K$. We have
firstly the following three lemmas, which establish the gain of
regularity in the velocity variable $v$, in the space variable $x$
and in the time variable $t$, respectively.

\begin{lem}\label{lemma2.0}
There exists a constant $C_{K}$ such that for any $f\in
C_0^\infty(K)$,
\begin{eqnarray*}
\|\nabla_vf\|_0\leq C_{K}\bigset{\abs{\inner{\mc{L}f,~f}}+\|f\|_0}.
\end{eqnarray*}
Moreover, for any $r\geq0,$ for any $\eps>0,$
\[\norm{\nabla_vf}_r\leq \eps\norm{\lo f}_r+C_{K,\eps}\norm{f}_r.\]

\end{lem}

We get a gain of regularity of order 1 for $v$ variable. This is
obtained directly by the positivity of coefficient $a$ and compact
support of $f$. For the space variable $x$, we have also the
following subelliptic estimate.

\begin{lem}\label{lemma2.1}
There exists a constant $C_K$ such that for any $f\in
C_0^\infty(K),$
\begin{eqnarray*}
\|D_x^{2/3}f\|_0\leq C_K(\|\mc{L}f\|_0+\|f\|_0),
\end{eqnarray*}
where $D_x^{2/3}=(-\triangle_x)^{1/3}$.
\end{lem}

This is a result of \cite{Bou}, and it is deduced by following  two
estimates
\[
\|D_x^{2/3}f\|_0\leq C_K~~\|\triangle_vf\|_0^{1/3}\|\partial_t
f+v\cdot\partial_x f\|_0^{2/3},
\]
and
\[
\|\triangle_vf\|_0\leq C_K(~~\|\mc L f\|_0+\|f\|_0~).
\]

For the time variable $t$, we have also a gain of regularity of
order $2/3$.

\begin{lem}\label{lemma2.2}
There exists a constant $C_K$ such that for any $f\in
C_0^\infty(K),$
\begin{eqnarray*}
\|\partial_tf\|_{-1/3}\leq C_K(\|\mc{L}f\|_0+\|f\|_0).
\end{eqnarray*}
\end{lem}
In fact, we have
$$
\|\partial_tf\|_{-1/3}=\|\Lambda^{-1/3}\partial_tf\|_0
\leq\|\Lambda^{-1/3}(\partial_t+v\cdot\partial_x)f\|_0
+\|\Lambda^{-1/3}v\cdot\partial_xf\|_0,
$$
where $\Lambda=(1+|D_t|^2+|D_x|^2+|D_v|^2)^{1/2}$. From Lemma
\ref{lemma2.1}, we have
$$
\|\Lambda^{-1/3}v\cdot\partial_xf\|_0\leq C_K\|D_x^{2/3}f\|_0 \leq
C_K(\|\mc{L}f\|_0+\|f\|_0).
$$
The estimation for the term $\|\Lambda^{-1/3}(\partial_t
+v\cdot\partial_x)f\|_0$ can be obtained by direct calculus as in
\cite{MX}.

\bigbreak \noindent {\bf{Proof of Proposition \ref{prp2.1}.}} The
Lemma \ref{lemma2.0}, Lemma \ref{lemma2.1} and Lemma \ref{lemma2.2}
deduce immediately
\begin{equation}
\label{Dx14}
\begin{array}{lll}
\|f\|_{2/3} &\leq& C_{K}\{~\|\mc{L}f\|_0+\|f\|_0~\}.
\end{array}
\end{equation}
Moreover, choose a function $\psi\in C_0^{\infty}(\bb R^{2n+1})$
such that $\psi|_K\equiv1$, Supp $\psi$ is a neighborhood of $K$.
Then  for any $f\in C_0^{\infty}( K)$ and  any $r\geq 0$,
\begin{eqnarray*}
\|f\|_r&=&\|\psi f\|_r\leq
C_K\{~\|\psi\Lambda^{r-{2/3}}f\|_{2/3}+\|[\Lambda^{r-{2/3}},~\psi]f\|_{2/3}~\}.
\end{eqnarray*}
By virtue of  (\ref{Dx14}) and the interpolation inequality
(\ref{interpolation}), we have
\begin{eqnarray*}
\|f\|_r&\leq&C_{K}\{~\|\mc{L}\psi\Lambda^{r-{2/3}}f\|_0+\|f\|_{r-{2/3}}~\}\\
&\leq&C_{\varepsilon,K}\{~\|\mc{L}\psi\Lambda^{r-{2/3}}f\|_0+\|f\|_0~\}+\varepsilon\|f\|_r.
\end{eqnarray*}
Taking $\varepsilon$ small enough, we get
\begin{eqnarray*}
\|f\|_r\leq C_K\{~\|\mc{L}f\|_{r-2/3}
+\|f\|_0+\|[\mc{L},~\psi\Lambda^{r-{2/3}}]f\|_0~\}.
\end{eqnarray*}
Direct verification gives
\begin{eqnarray*}
[\mc{L},~\psi\Lambda^{r-{2/3}}]
&=&[\partial_t+v\cdot\partial_x,~\psi\Lambda^{r-{2/3}}]-\sum\limits_{j=1}^n\{~\|[a,~\psi\Lambda^{r-{2/3}}]\partial_{v_j}^2\\
&&+a[\partial_{v_j},~[\partial_{v_j},\psi\Lambda^{r-{2/3}}]~]
+2a[\partial_{v_j},~\psi\Lambda^{r-{2/3}}]\partial_{v_j}~\},
\end{eqnarray*}
This along with Lemma \ref{lemma2.0}  yields
\begin{eqnarray*}
\|[\mc{L},~\psi\Lambda^{r-{2/3}}]f\|_0&\leq& C_K\{~\|f\|_{r-{2/3}}
+\sum\limits_{j=1}^n \|\partial_{v_j}f\|_{r-{2/3}}~~\}\\
&\leq&C_K\{~\|\mc Lf\|_{r-2/3}+\|f\|_{r-2/3}~\}.
\end{eqnarray*}
These three estimates gives immediately
\begin{eqnarray*}
\|f\|_r\leq C_K\{~\|\mc{L}f\|_{r-2/3}+\|f\|_0+\|f\|_{r-2/3}~\}.
\end{eqnarray*}
Applying interpolation inequality (\ref{interpolation}) again and
taking $\varepsilon$ small
enough, we prove Proposition \ref{prp2.1}.\\

We consider now the commutators of the operators $\mc{L}$ with
derivation and cut-off function.

\begin{prp}\label{prp2.2}
Let $K$ be a compact subset of $\bb{R}^{2n+1}.$  Then for any $r\geq
0,$ there exist constants $C_{K,r}, C_{K,r,\varphi}$ such that for
any $f\in C^\infty_0(K),$
\begin{eqnarray*}
\|[\mc L, ~D]f\|_r
&\leq&C_{K,r}\{~\|\mc{L}f\|_{r+1-{2/3}}+\|f\|_0~\},
\end{eqnarray*}
and
\begin{equation*}
\|[\mc L,~\varphi]f\|_r \leq C_{K,r,\varphi}\{~\|\mc
Lf\|_{r-1/3}+\|f\|_0~\},
\end{equation*}
where $\varphi\in C^\infty_b(\bb R^{2n+1})$ and we denote by $D$ the
differential operator $\partial_t,\partial_x$ or $\partial_v.$
\end{prp}

\noindent {\textbf{Proof.} \hspace{0.1cm} By using the positivity of
coefficient $a$, we have
\begin{eqnarray*}
\|\triangle_vf\|_r &\leq& C_{K}\{~\|\mc Lf\|_r+\|f\|_{r+1} ~\}.
\end{eqnarray*}
And $[\mc L, ~D]=[\partial_t+v\cdot\partial_x,~D]-[a,
~D]\triangle_v$ deduce
\begin{eqnarray*}
\|[\mc L, ~D]f\|_r \leq C_{K}\{~\|f\|_{r+1}+\|\triangle_vf\|_r~\}.
\end{eqnarray*}
The above two inequalities  along with the subelliptic estimate
(\ref{subestimate}) yield the first desired inequality in
Proposition \ref{prp2.2}.

To treat $\|[\mc L,~\varphi]f\|_r$, the subelliptic estimate
(\ref{subestimate}) give
\begin{eqnarray*}
\|\nabla_vf\|_r \leq C_K (\|\mc Lf\|_{r-1/3}+\|f\|_{0}).
\end{eqnarray*}
Now simple verification gives
\begin{eqnarray*}
\|[\mc L,~\varphi]f\|_r&\leq&C_K\big\{~\|f\|_r+\sum\limits_{j=1}^n~\|\partial_{v_j}f\|_r~\big\}\\
&\leq& C_{K,r} \big\{~\|\mc Lf\|_{r-1/3}+\|f\|_{0}~\big\}.
\end{eqnarray*}
This completes the proof of  Proposition \ref{prp2.2}.

\bigskip We prove now the Gevrey hypoellipticity of $\mc{L}$ . Our
starting point is the following result due to M.Durand [\ref{Dur}]:

\begin{prp}\label{prp1}
Let $P$ be a linear differential operator with smooth coefficients
in $\bb{R}_{y}^m$ and $\varrho ,\varsigma $  two fixed positive
numbers. If for any $r \geq 0$, any compact $K\subseteq {\bb{R}^m}$
and any $\varphi\in C^\infty(\bb{R}^m)$, there exist constants
$C_{K,r}$ and $C_{K,r}(\varphi)$ such that for all $f\in
C_0^\infty(K)$,the following conditions are fulfilled:
\begin{eqnarray*}
(H_1) \hspace{5.5cm}\|f\|_r\leq
C_{K,r}(\|Pf\|_{r-\varrho}+\|f\|_0),\hspace{6.1cm}\\
(H_2) \hspace{4cm}\|[P,~D_j]f\|_r\leq
C_{K,r}(\|Pf\|_{r+1-\varsigma}+\|f\|_0),\hspace{6cm}\\
(H_3) \hspace{4.1cm}\|[P,~\varphi]f\|_r\leq
C_{K,r}(\varphi)(\|Pf\|_{r-\varsigma}+\|f\|_0),\hspace{6cm}
\end{eqnarray*}
where$$D_j=\frac{1}{i}\frac{\partial}{\partial
y_j},j=1,2,\cdots,m.$$Then for $s\geq \max(1/\varsigma,2/\varrho),P$
is $G^s(\bb{R}^m)$ hypoelliptic, provided the coefficients of $P$
are in the class of $G^s(\bb{R}^m).$
\end{prp}

Proposition \ref{prp2.1} shows that the operator $\mc L$ satisfies
the conditions $(H_1)$ with $\varrho=2/3$, Proposition \ref{prp2.2}
assures the conditions $(H_2)$ and $(H_3)$ with $\varsigma=1/3$.
Then $\mc L$ is $G^s(\bb{R}^{2n+1})$ hypoelliptic, $s\geq 3$, and we
have proved Theorem \ref{th1}.


\section{Gevrey regularity of nonlinear equations}
\label{sect3} \setcounter{equation}{0}

Let $u\in L^{\infty}_{loc}(\bb R^{2n+1})$ be a weak solution of
(\ref{++1.2}). Firstly, we will prove $u\in C^\infty(\bb R^{2n+1}).$
And we need the following nonlinear composition results (see for
example \cite{Xu}).

\begin{lem}\label{composition}

Let $F(t,x,v,w,p)\in C^\infty(\bb R^{2n+2+n})$ and $r\geq 0$. If $u,
\nabla_vu\in L^{\infty}_{loc}(\bb R^{2n+1})\cap H_{loc}^{r}(\bb
R^{2n+1})$, then $F\big(\cdot,u(\cdot),\nabla_vu(\cdot)\big)\in
H_{loc}^{r}(\bb R^{2n+1}),$ and
\begin{equation}\label{F:regularity}
 \norm{\phi_1F\big(\cdot,u(\cdot),\nabla_vu(\cdot)\big)}_r
 \leq \bar C \set{~\norm{\phi_2u}_r+\norm{\phi_2 \nabla_vu}_r~},
\end{equation}
where $\phi_1, \phi_2 \in C_0^{\infty}(\bb R^{2n+1})$ and $\phi_2=1
$ on the support of $\phi_1$, and $\bar C$ is a constant depending
only on $r, \phi_1,\phi_2.$

\end{lem}

\noindent{\bf Remark.} If the nonlinear term $F$ is independent of
$p$ or in the form of
$$
\nabla_v (F(t, x, v, u)),
$$
Then that $u \in L^{\infty}_{loc}(\bb R^{2n+1})\cap H_{loc}^{r}(\bb
R^{2n+1})$ yields $F\big(\cdot,u(\cdot),\nabla_vu(\cdot)\big)\in
H_{loc}^{r}(\bb R^{2n+1}).$

\begin{lem}\label{priori estimate}

Let $u,\nabla_v u\in H_{loc}^r(\Real^{2n+1}), r\geq0$. Then we have
\begin{equation}\label{priori}
 \norm{\vpi_1\nabla_v u}_r\leq C \norm{\vpi_2u}_r,
\end{equation}
where $\vpi_1, \vpi_2\in C_0^\infty(\Real^{2n+1})$ and  $\vpi_2=1$
on the support of  $\vpi_1$,  and $C$ is a constant depending only
on $r, \vpi_1,\vpi_2.$

\end{lem}

In fact, we have
 \[\norm{\vpi_1\nabla_vu}_r
 \leq\norm{[\nabla_v,~\vpi_1] u}_r+\norm{\nabla_v\vpi_1 u}_r.\]
Clearly, the first term on the right is bounded by
$C\norm{\vpi_2u}_r$. For the second term , combining the second
inequality in Lemma \ref{lemma2.0} and \reff{F:regularity}, we get
the desired estimate \reff{priori} at once. This completes the proof
of Lemma \ref{priori estimate}.

\bigskip Now we are ready to prove

\begin{prp}\label{+smooth}

Let $u$ be a weak solution of (\ref{++1.2}) such that $u,\nabla_v
u\in L_{loc}^\infty(\Real^{2n+1})$. Then $u$ is in $C^\infty(\bb
R^{2n+1})$.

\end{prp}

In fact, from the subelliptic estimate (\ref{subestimate}) and the
fact $\mc{L}u(\cdot)=F(\cdot,u(\cdot),\nabla_vu(\cdot))$, it then
follows
\begin{eqnarray}\label{+estimate}
\|\psi_1u\|_{r+2/3}\leq \bar C \{~\|\psi_2
F\big(\cdot,u(\cdot),\nabla_vu(\cdot)\big)\|_r+\|\psi_2 u\|_0~\},
\end{eqnarray}
where $\psi_1, \psi_2 \in C_0^{\infty}(\bb R^{2n+1})$ and $\psi_2=1
$ on the support of $\psi_1$.  Combining \reff{F:regularity},
\reff{priori} and \reff{+estimate}, we have $u\in
H_{loc}^{\infty}(\bb R^{2n+1})$ by standard iteration. This
completes the proof of Proposition \ref{+smooth}.

Now starting from the smooth solution, we prove the Gevrey
regularity. It suffices to show the  regularity in the open unit
ball
$$\Omega=\{(t,x,v)\in \bb{R}^{2n+1}: t^2+|x|^2+|v|^2<1 \}.$$
Set
\[\Omega_\rho=\big\{(t,x,v)\in \Omega: \big(t^2+|x|^2+|v|^2\big)^{1/2}<1-\rho \big\},
\hspace{0.6cm}0<\rho<1.\]

Let $U$ be an open subset of $\bb{R}^{2n+1}$. Denote by $H^r(U)$ the
space consisting of the functions which are defined in $U$ and can
be extended to $H^r(\bb{R}^{2n+1})$. Define
$$\|u\|_{H^r(U)}=\inf\big\{\|\tilde{u}\|_{H^s(\mathbb{R}^{n+1})}:\tilde{u}\in
H^s(\mathbb{R}^{2n+1}),\tilde{u}|_U=u\big\}.$$ We denote
$\|u\|_{r,U}=\|u\|_{H^r(U)},$ and
$$\|D^ju\|_r=\sum_{|\beta|=j}\|D^\beta u\|_r.$$

In order to treat the nonlinear term $F$ on the right hand of
(\ref{++1.2}), we need the following two lemmas. The first one (see
\cite{Xu} for example) concerns weak solution in some algebra, and
the second  is an analogue of Lemma 1 in \cite{Friedman}. In the
sequel $C_j>1$ will be used to denote suitable constants depending
only on $n$ or the function $F$.

\begin{lem}\label{algebra}

Let $r>(2n+1)/2$ and  $u_1, u_2\in H^{r}(\bb R^{2n+1})$, Then
$u_1u_2\in H^r(\bb R^{2n+1})$, moreover
\begin{eqnarray}\label{4.1}
\|u_1u_2\|_r\leq \tilde C \|u_1\|_r\|u_2\|_r,
\end{eqnarray}
where $\tilde C$ is a constant depending only on $n,r.$

\end{lem}

\begin{lem}\label{+Friedman}

Let $M_j$ be a sequence of positive numbers and for some $B_0>0,$
the $M_j$ satisfy the monotonicity conditions
\begin{equation}\label{monotonicity}
\frac{j!}{i!(j-i)!}M_iM_{j-i}\leq B_0M_j, ~~(i=1,2,\cdots,j;~~
j=1,2,\cdots).
\end{equation}
Suppose $F(t,x,v,u,p)$ satisfy
\begin{equation}\label{condition1}
\big\|\left(D_{t,x,v}^{j}D_u^lD_p^{m}F\right)\big(\cdot, u(\cdot),
\nabla_vu(\cdot)\big)\big\|_{r+n+1,\Omega}\leq C_1^{j+l+m}
M_{j-2}M_{m+l-2}, \quad j,m+l\geq2,
\end{equation}
where $r$ is a real number satisfying $r+n+1>(2n+1)/2$. Then there
exist two constants $C_2, C_3$ such that for any $H_0,H_1$
satisfying $H_0, H_1\geq1$ and $H_1\geq
 C_2 H_0$, if $u(t,x,v)$ satisfy the following conditions
\begin{equation}\label{condition2}
  \|D^j u\|_{r+n+1, \Omega_{\tilde\rho}}
  \leq H_0, \quad 0\leq j\leq1,
\end{equation}
\begin{equation}\label{condition3}
  \|D^j u\|_{r+n+1,\Omega_{\tilde\rho}}
  \leq H_0H_1^{j-2}M_{j-2},\quad 2\leq j\leq N,
\end{equation}
\begin{equation}\label{condition1}
  \|D_vD^j u\|_{r+n+1,\Omega_{\tilde\rho}}\leq
  H_0H_1^{j-2}M_{j-2},\quad 2\leq j\leq N.
\end{equation}
 Then for all $\alpha$ with $|\alpha|=N$,
\begin{eqnarray}\label{conclusion}
\big\|\psi_{N}D^{\alpha}\big[F\big(\cdot,
u(\cdot),\nabla_vu(\cdot)\big)\big]\big\|_{r+n+1}\leq C_3
H_0H_1^{N-2}M_{N-2},
\end{eqnarray}
where $\psi_N\in C_0^\infty(\Omega_{\tilde\rho})$ is an arbitrary
function.

\end{lem}

\noindent{\bf Proof.} Denote $p=(p_1, p_2,\cdots, p_n)=\nabla_vu$
and $k=(k_1, k_2\cdots, k_n).$  From Faa di Bruno' formula,
$\psi_{N}D^\alpha[F(\cdot,u(\cdot),\nabla_v u(\cdot))]$ is the
linear combination of terms of the form
\begin{equation}\label{++A}
 \frac{\psi_{N}\partial^{|\tilde\alpha|+l+|k|}F}{\partial_{t,x,v}^{\tilde\alpha}
 \partial u^l\partial p_1^{k_1}\cdots\partial p_n^{k_n}}
 \prod_{j=1}^lD^{\gamma_j}u\cdot\prod_{i=1}^n~\prod_{j_i=1}^{k_i}D^{\beta_{j_i}}
 (\partial_{v_i}u),
\end{equation}
where $|\tilde\alpha|+l+|k|\leq |\alpha|$ and
\[
\sum_{j=1}^l\gamma_i+\sum_{i=1}^n\sum_{j_i}^{k_i}\beta_{j_i}=\alpha-\tilde\alpha,\]
and if $\gamma_i$ or $\beta_{j_i}$ equals to 0, we just mean
$D^{\gamma_i}u$ or $D^{\beta_{j_i}}u$ doesn't appear in (\ref{++A}).
Choose a function $\tilde\psi\in C_0^\infty(\Omega_{\tilde\rho})$
such that $\tilde\psi=1$ on Supp $\psi_N$. Note that
$n+1+r>(2n+1)/2,$ and hence applying Lemma \ref{algebra}, we have
\begin{equation}\label{++a}
\begin{array}{lll}
  &&\norm{\frac{\psi_{N}\partial^{|\tilde\alpha|+l+|k|}F}{\partial_{t,x,v}^{\tilde\alpha}
    \partial u^l\partial p_1^{k_1}\cdots\partial p_n^{k_n}}\prod_{j=1}^lD^{\gamma_j}u
    \cdot\prod_{i=1}^n~\prod_{j_i=1}^{k_i}D^{\beta_{j_i}}(\partial_{v_i}u)}_{r+n+1}\\\\
  &=&\norm{\frac{\psi_{N}\partial^{|\tilde\alpha|+l+|k|}F}{\partial_{t,x,v}^{\tilde\alpha}
    \partial u^l\partial p_1^{k_1}\cdots\partial p_n^{k_n}}\prod_{j=1}^l\tilde\psi D^{\gamma_j}u
    \cdot\prod_{i=1}^n~\prod_{j_i=1}^{k_i}\tilde\psi\partial_{v_i}D^{\beta_{j_i}}u}_{r+n+1}\\\\
  &\leq&  \tilde C\norm{\psi_N(\partial^{|\tilde\alpha|+l+|k|}F)}_{r+n+1}\cdot
     \prod_{j=1}^{l}\norm {\tilde\psi D^{\gamma_j}u}_{r+n+1}\times
     \prod_{i=1}^n~\prod_{j_i=1}^{k_i}\norm{\tilde\psi\partial_{v_i}D^{\beta_{j_i}}u}_{r+n+1}\\\\
  &\leq&  C_0\norm{(\partial^{|\tilde\alpha|+l+|k|}F)}_{r+n+1,\Omega}\cdot
     \prod_{j=1}^{l}\norm{ D^{\gamma_j}u}_{ r+n+1,\Omega_{\tilde\rho}}\times
     \prod_{i=1}^n~\prod_{j_i=1}^{k_i}\norm{\partial_{v_i}D^{\beta_{j_i}}u}
     _{r+n+1,\Omega_{\tilde\rho}}.
\end{array}
\end{equation}
In virtue of (\ref{condition2})-(\ref{condition1}) and (\ref{++a}),
the situation is entirely similar to \cite{Friedman}. The only
difference is that we replace the H\"{o}lder norm  $|u|_j$  by
$\|D^ju\|_{ r+n+1, \Omega_{\tilde\rho}}$ and
$\norm{D_vD^ju}_{r+n+1,\Omega_{\tilde\rho}}$. Then the same argument
as the proof of Lemma 1 in \cite{Friedman} yields
(\ref{conclusion}). This completes the proof of Lemma
\ref{+Friedman}.

\smallskip
\begin{prp}\label{prp4'}

Let $s\geq3$.  Suppose $u\in C^\infty(\bar\Omega)$ is a solution of
(\ref{++1.2}), and $a(t,x,v)\in G^s(\mathbb{R}^{2n+1})$, $ F(t, x,
v, w,p)\in G^s(\mathbb{R}^{2n+2+n})$ and $a\geq c_0>0$. Then there
exits a constant $A$ such that for any $r\in[0,1]$ and any
$N\in\bb{N}$, $N\geq3,$
\begin{eqnarray*}\label{rk2}
(E)_{r, N}\hspace{0.6cm}&&\|D^\alpha u\|_{r+n+1,\Omega_\rho}+\|D_v
D^\alpha u\|_{r-1/3+n+1,\Omega_{\rho}}\\ &\leq&
\frac{A^{|\alpha|-1}}{\rho^{s(|\alpha|-3)}}\big((|\alpha|-3)!\big)^{s}(N/\rho)^{sr},
\quad \forall ~|\alpha|= N,~~\forall~ 0<\rho< 1.
\end{eqnarray*}

\end{prp}

\bigbreak  From $(E)_{r,N}$ , we have immediately

\begin{prp}\label{prp4}

Under the same assumption as Proposition \ref{prp4'},  we have $u\in
G^s({\Omega}).$
\end{prp}

In fact, for any compact sunset $K$ of $\Omega$, we have
$K\subset\Omega_{\rho_0}$ for some $\rho_0, ~0<\rho_0<1$. For any
$\alpha, ~~|\alpha|\geq3, $ letting $r=0$ in $(E)_{r,N}$, we have
\begin{equation*}
\begin{array}{lll}\|D^\alpha u\|_{L^2(K)}&\leq&\|D^\alpha
u\|_{n+1,\Omega_{\rho_0}} \leq
\frac{A^{|\alpha|-1}}{{\rho_0}^{s(|\alpha|-3)}}\big((|\alpha|-3)!\big)^{s}
\leq \big({A\over{\rho_0}^s}\big)^{|\alpha|+1}(|\alpha|!)^s.
\end{array}
\end{equation*}
This completes the proof of Proposition \ref{prp4}.

\bigbreak \noindent{\bf Proof of Proposition \ref{prp4'}.}  We use
induction on $N$. Assuming $(E)_{r,N-1}$ holds for any $r$ with
$0\leq r\leq1$, and we will show $(E)_{r,N}$ still holds for any
$r\in[0,~1]$. For any $\alpha, |\alpha|= N,$ and for any
$\rho,~0<\rho<1,$ choose a function $\varphi_{\rho,N}\in
C_0^\infty(\Omega_{{{(N-1)\rho}\over N}})$ such that
$\varphi_{\rho,N}=1$ in $\Omega_{\rho}$.  it is easy to see
  \[  \sup|D^\gamma\varphi_{\rho,N}|\leq C_\gamma (\rho/N)^{-|\gamma|}\leq
      C_\gamma (N/\rho)^{|\gamma|},\indent \forall~\gamma.
  \]
 And we will proceed to prove the truth of $(E)_{r,N}$
by the following lemmas.

\begin{lem}\label{r0}

For $r=0$,  we have
\[\|D^\alpha u\|_{n+1, \Omega_\rho}+\|D_v
D^\alpha u\|_{-1/3+n+1,  \Omega_\rho}\leq
\frac{C_7A^{|\alpha|-2}}{\rho^{s(|\alpha|-3)}}\big((|\alpha|-3)!\big)^{s},
\quad \forall~~0<\rho<1.\]

\end{lem}

\noindent{\bf Proof .} Write $|\alpha|=|\beta|+1$, then
$|\beta|=N-1$. Denote ${{N-1}\over N}\rho$ by $\tilde \rho.$ In the
sequel we will use the following fact frequently
\[{1\over{\rho}^{sk}}\leq{1\over{\tilde\rho}^{sk}}={1\over{\rho}^{sk}}\times\big({N\over{N-1}}\big)^{sk}
\leq  {C_4\over{\rho}^{sk}},\quad k=1,2,\cdots, N-3.\] Note that
$\varphi_{\rho,N}=1$ in $\Omega_\rho$ and hence
\begin{eqnarray*}
\|D^\alpha u\|_{n+1,\Omega_\rho}&\leq&\|\varphi_{\rho,N}D^\alpha
u\|_{n+1} \leq\|\varphi_{\rho,N}D^\beta
u\|_{1+n+1}+\|(D\varphi_{\rho,N})D^\beta u\|_{n+1}\\&\leq&
C_5\{~\|D^\beta u\|_{1+n+1,\Omega_{\tilde \rho}}+(N/\rho)\|D^\beta
u\|_{n+1,\Omega_{\tilde \rho}}~\}.
\end{eqnarray*}
Since $(E)_{r,N-1}$ holds by assumption for any $r$ with $0\leq
r\leq1$ , we have immediately
\begin{eqnarray*}
  &&\|D^\beta u\|_{1+n+1,\Omega_{\tilde\rho}}
     +(N/\rho)\|D^\beta u\|_{n+1,\Omega_{\tilde\rho}}\\
  &\leq&\frac{A^{|\beta|-1}}{\tilde\rho^{s(|\beta|-3)}}
     \big((|\beta|-3)!\big)^{s}(N/\tilde\rho)^{s}
     +(N/\rho)\frac{A^{|\beta|-1}}{\tilde\rho^{s(|\beta|-3)}}\big((|\beta|-3)!\big)^{s}\\
  &\leq&\frac{2A^{|\alpha|-2}}{{\tilde\rho}^{s(|\alpha|-3)}}
     \big((|\alpha|-3)!\big)^{s}\big(N/(N-3)\big)^{s}\\
  &\leq&\frac{C_6A^{|\alpha|-2}}{\rho^{s(|\alpha|-3)}}\big((|\alpha|-3)!\big)^{s}.
\end{eqnarray*}
Thus
\begin{eqnarray}\label{+inequ}
\|D^\alpha u\|_{n+1, \Omega_\rho}\leq
\frac{C_5C_6A^{|\alpha|-2}}{\rho^{s(|\alpha|-3)}}\big((|\alpha|-3)!\big)^{s}.
\end{eqnarray}
The same arguments as above shows that
\begin{eqnarray*}
\|D_v D^\alpha u\|_{-1/3+n+1, \Omega_\rho} \leq
\frac{C_5C_6A^{|\alpha|-2}}{\rho^{s(|\alpha|-3)}}\big((|\alpha|-3)!\big)^{s}.
\end{eqnarray*}
This along with (\ref{+inequ}) yields the conclusion.

\begin{lem}\label{r1/3}

For any $0\leq r\leq 1/3$, we have
\[\|D^\alpha u\|_{r+n+1, \Omega_\rho}+\|D_v
D^\alpha u\|_{r-1/3+n+1,  \Omega_\rho}\leq
\frac{C_{35}A^{|\alpha|-2}}{\rho^{s(|\alpha|-3)}}\big((|\alpha|-3)!\big)^{s}(N/\rho)^{rs},
\hspace{0.3cm} \forall~~0<\rho<1.\]

\end{lem}

\noindent{\bf Proof .} We firstly prove the conclusion is true for
 $ r=1/3$, i.e., to show
\begin{equation*}\label{r1/3'}
\begin{array}{rrr}
  \|D^\alpha u\|_{1/3+n+1, \Omega_\rho}+\|D_vD^\alpha u\|_{n+1, \Omega_\rho}
  \leq\frac{C_{35}A^{|\alpha|-2}}{\rho^{s(|\alpha|-3)}}\big((|\alpha|-3)!\big)^{s}
     (N/\rho)^{s/3},\quad \forall~~0<\rho<1.
\end{array}
\end{equation*}
And we will proceed in the following four steps.

\bigbreak {\bf \emph{Step 1.}} Claim
\begin{equation}
\label{step1}
\begin{array}{lll}
\|[\mc{L},~~ \varphi_{\rho,N}D^\alpha] u\|_{-1/3+n+1}\leq
\frac{C_{19}A^{|\alpha|-2}}{\rho^{s(|\alpha|-3)}}\big((|\alpha|-3)!\big)^{s}(N/\rho)^{s/3}.
\end{array}
\end{equation}

In fact, write $\mc L=X_0-a\triangle_v$ with
$X_0=\partial_t+v\cdot\partial_x$. Then direct verification deduces
\begin{eqnarray*}
\|[\mc{L},~~ \varphi_{\rho,N}D^\alpha] u\|_{-1/3+
n+1}&\leq&\|[X_0,~~ \varphi_{\rho,N}D^\alpha]
u\|_{-1/3+n+1}+\|a[\triangle_v,~~ \varphi_{\rho,N}D^\alpha]
u\|_{-1/3+ n+1}\\&&+\|\varphi_{\rho,N}[a,~~
D^\alpha]\triangle_v u\|_{-1/3+n+1}\\
&=:&(I)+(II)+(III).
\end{eqnarray*}
Denote $[X_0,~~ D^\alpha]$ by $D^{\alpha_0}$. Then $|\alpha_0|\leq
|\alpha|$ and
\begin{eqnarray*}
(I) &\leq&\|[X_0,~~ \varphi_{\rho,N}]D^\alpha
u\|_{n+1}+\|\varphi_{\rho,N} D^{\alpha_0}
u\|_{n+1}\\
&\leq&C_8\big\{~(N/\rho)\|D^\alpha u\|_{n+1,\Omega_{\tilde
\rho}}+\|D^{\alpha_0} u\|_{n+1,\Omega_{\tilde \rho}}~\big\}.
\end{eqnarray*}
Note that $s\geq3$. Using Lemma \ref{r0}, we have
\begin{equation}\label{+I}
\begin{array}{lll} (I) &\leq&
C_8\big(N/\rho+1\big)\frac{C_7A^{|\alpha|-2}}{{\tilde\rho}^{s(|\alpha|-3)}}\big((|\alpha|-3)!\big)^{s}\leq
\frac{C_9A^{|\alpha|-2}}{\rho^{s(|\alpha|-3)}}\big((|\alpha|-3)!\big)^{s}(N/\rho)^{s/3}.
\end{array}
\end{equation}
Next we will estimate $(II)$. It is easy to see that
\begin{equation}\label{IIa}
\begin{array}{lll}
\|[\triangle_v,~\varphi_{\rho,N}]D^\alpha u\|_{-1/3+n+1}&\leq&
2\|[D_v,~\varphi_{\rho,N}]D_vD^\alpha
u\|_{-1/3+n+1}\\\\
&&+\|[D_v,~[D_v,~\varphi_{\rho,N}]~]D^\alpha u\|_{-1/3+ n+1}.
\end{array}
\end{equation}
We firstly treat the first term on the right hand. Using Lemma
\ref{r0} again, we have
\begin{align}\label{IIb}
\begin{split}
   \|[D_v,~\varphi_{\rho,N}]D_vD^\alpha u\|_{-1/3+n+1}
   &\leq(N/\rho)\|D_vD^\alpha u\|_{-1/3+n+1,
   \Omega_{\tilde\rho}}\\
   &\leq(N/\rho)\frac{C_7A^{|\alpha|-2}}{{\tilde\rho}^{s(|\alpha|-3)}}
     \big((|\alpha|-3)!\big)^{s}\\
   &\leq\frac{C_{10}A^{|\alpha|-2}}{{\rho}^{s(|\alpha|-3)}}
     \big((|\alpha|-3)!\big)^{s}(N/\rho)^{s/3}.
\end{split}
\end{align}
Next we treat $\|[D_v,~[D_v,~\varphi_{\rho,N}]~]D^\alpha
u\|_{-1/3+n+1}$, and we compute
\begin{eqnarray*}
  &&\|[D_v,~[D_v,~\varphi_{\rho,N}]~]D^\alpha u\|_{-1/3+n+1}\\
  &\leq&\|(D^2\varphi_{\rho,N})D^\beta u\|_{2/3+n+1}
    +\|(D^3\varphi_{\rho,N})D^\beta u\|_{-1/3+n+1}\\
  &\leq&C_{11}\big\{~(N/\rho)^2\|D^\beta u\|_{2/3+n+1,\Omega_{\tilde\rho}}
    +(N/\rho)^3\|D^\beta u\|_{n+1,\Omega_{\tilde\rho}}~\big\}\\
  &\leq&C_{11}\big\{~(N/\rho)^2{A^{|\beta|-1}\over{\tilde\rho^{s(|\beta|-3)}}}
    \big((|\beta|-3)!\big)^{s}(N/\tilde\rho)^{2s/3}\\
    &&+(N/\rho)^3{A^{|\beta|-1}\over{\tilde\rho^{s(|\beta|-3)}}}
    \big((|\beta|-3)!\big)^{s}~\big\}\\
  &\leq&C_{11}\big\{~(N/\rho)^2(N/\tilde\rho)^{-s/3}
    {A^{|\alpha|-2}\over{\tilde\rho^{s(|\alpha|-3)}}}\big((|\alpha|-3)!\big)^{s}\\
    &&+(N/\rho)^3(N/\tilde\rho)^{-s}{A^{|\alpha|-2}\over{\tilde\rho^{s(|\alpha|-3)}}}
    \big((|\alpha|-3)!\big)^{s}~\big\}\\
 &\leq&{C_{12}A^{|\alpha|-2}\over{\rho^{s(|\alpha|-3)}}}\big((|\alpha|-3)!\big)^{s}
   (N/\rho)^{s/3}.
\end{eqnarray*}
This along with (\ref{IIa}) and (\ref{IIb}) shows at once
\begin{eqnarray}\label{+II}
(II)
\leq\frac{C_{13}A^{|\alpha|-2}}{\rho^{s(|\alpha|-3)}}\big((|\alpha|-3)!\big)^{s}(N/\rho)^{s/3}.
\end{eqnarray}
It remains to treat $(III)$, and using Leibniz' formula,
\begin{eqnarray*}
  (III)&\leq& \sum\limits_{0<|\gamma|\leq|\alpha|}
        \big(\begin{array}{c}\alpha\\ \gamma\end{array}\big)
        \big\|\varphi_{\rho,N}(D^{\gamma}a)\triangle_vD^{\alpha-\gamma} u\big\|_{-1/3+n+1}\\
  &\leq& \sum\limits_{0<|\gamma|\leq|\alpha|}
       \big(\begin{array}{c}\alpha\\ \gamma\end{array}\big)
       \big\|D^{\gamma}a\|_{n+1,\Omega}\cdot
       \|\varphi_{\rho,N}\triangle_vD^{\alpha-\gamma} u\big\|_{-1/3+n+1}.
\end{eqnarray*}
Since $a\in G^s(\bb R^{2n+1})$, then
   \[ \|D^{\gamma}a\|_{n+1, \Omega}
     \leq C_{14}^{|\gamma|-2}\big((|\gamma|-3)!\big)^{s},
     \quad|\gamma|\geq3,\]
and
   \[\|D^\gamma a\|_{n+1, \Omega}
     \leq C_{14},\quad|\gamma|=1,~2.\]
Moreover, note $|\alpha|-|\gamma|+1\leq N,$ and hence applying Lemma
\ref{r0}, we have for any $\gamma, ~|\gamma|\leq|\alpha|-2$,
\begin{eqnarray*}
  \|\varphi_{\rho,N} \triangle_vD^{\alpha-\gamma} u\|_{-1/3+n+1}
  &\leq&\|D_vD^{\alpha-\gamma+1}u\|_{-1/3+n+1,\Omega_{\tilde\rho}}\\
  &\leq&\frac{C_{7}A^{|\alpha|-|\gamma|+1-2}}{\tilde\rho^{s(|\alpha|-|\gamma|-2)}}
    \big((|\alpha|-|\gamma|-2)!\big)^{s}\\
  &\leq&\frac{C_{15}A^{|\alpha|-|\gamma|+1-2}}{\rho^{s(|\alpha|-|\gamma|-2)}}
    \big((|\alpha|-|\gamma|-2)!\big)^{s}.
\end{eqnarray*}
Consequently, we compute
\begin{eqnarray*}
   &&\sum\limits_{2\leq|\gamma|\leq|\alpha|-2}
      \big(\begin{array}{c}\alpha\\ \gamma\end{array}\big)
      \big\|D^{\gamma}a\|_{n+1,\Omega}\cdot\|\varphi_{\rho,N}
      \triangle_vD^{\alpha-\gamma} u\big\|_{-1/3+n+1}\\
   &\leq&\sum\limits_{2\leq|\gamma|\leq|\alpha|-2}
     \big(\begin{array}{c}\alpha\\ \gamma\end{array}\big)
     C_{14}^{|\gamma|-2}\big((|\gamma|-2)!\big)^{s}
     \frac{C_{15}A^{|\alpha|-|\gamma|+1-2}}{\rho^{s(|\alpha|-|\gamma|-2)}}
     \big((|\alpha|-|\gamma|-2)!\big)^{s}\\
  &\leq&\frac{C_{15}A^{|\alpha|-2}}{\rho^{s(|\alpha|-3)}}
     \sum\limits_{2\leq|\gamma|\leq|\alpha|-2}
     \big({C_{14}\over A}\big)^{|\gamma|-1}|\alpha|!
     \big((|\gamma|-2)!\big)^{s-1}\big((|\alpha|-|\gamma|-2)!\big)^{s-1}\\
  &\leq&\frac{C_{15}A^{|\alpha|-2}}{\rho^{s(|\alpha|-3)}}\big((|\alpha|-3)!\big)^{s}
     \sum\limits_{2\leq|\gamma|\leq|\alpha|-2}
     \big({C_{14}\over A}\big)^{|\gamma|-1}|\alpha|
     \frac{(|\alpha|-1)(|\alpha|-2)}{(|\alpha|-3)^{s-1}}\\
   &\leq&\frac{C_{16}A^{|\alpha|-2}}{\rho^{s(|\alpha|-3)}}
     \big((|\alpha|-3)!\big)^{s}(N/\rho)^{s/3}
     \sum\limits_{2\leq|\gamma|\leq|\alpha|-2}
     \big({C_{14}\over A}\big)^{|\gamma|-1}.
\end{eqnarray*}
Taking $A$ large enough such that
$\sum\limits_{2\leq|\gamma|\leq|\alpha|-2} \big({C_{14}\over
A}\big)^{|\gamma|-1}\leq1$,
then we get
    \[ \sum\limits_{2\leq|\gamma|\leq|\alpha|-2}
       \big(\begin{array}{c}\alpha\\ \gamma\end{array}\big)
       \big\|D^{\gamma}a\|_{n+1,\Omega}\cdot\|\varphi_{\rho,N}
       \triangle_vD^{\alpha-\gamma} u\big\|_{-1/3+n+1}
       \leq\frac{C_{16}A^{|\alpha|-2}}{\rho^{s(|\alpha|-3)}}
       \big((|\alpha|-3)!\big)^{s}(N/\rho)^{s/3}.
    \]
For $|\gamma|=1, ~|\alpha|-1$ or $|\alpha|$, we can compute directly
    \[ \big(\begin{array}{c}\alpha\\ \gamma\end{array}\big)
       \big\|D^{\gamma}a\|_{n+1,\Omega}\cdot\|\varphi_{\rho,N}
       \triangle_vD^{\alpha-\gamma} u\big\|_{-1/3+n+1}
       \leq\frac{C_{17}A^{|\alpha|-2}}{\rho^{s(|\alpha|-3)}}
       \big((|\alpha|-3)!\big)^{s}(N/\rho)^{s/3}.
    \]
Combination of the above two inequalities give that
  \[ (III)\leq \frac{C_{18}A^{|\alpha|-2}}{\rho^{s(|\alpha|-3)}}
     \big((|\alpha|-3)!\big)^{s}(N/\rho)^{s/3}.
  \]
This along with (\ref{+I}) and (\ref{+II}) yields the conclusion
(\ref{step1}).

\bigbreak {\bf \emph{Step 2.}} Claim
\begin{equation}\label{step2}
\begin{array}{lll}
\|
\varphi_{\rho,N}D^\alpha[F\big(\cdot,u(\cdot),\nabla_vu(\cdot)\big)]\|_{
-1/3+n+1}\leq
\frac{C_{21}A^{|\alpha|-2}}{\rho^{s(|\alpha|-3)}}\big((|\alpha|-3)!\big)^{s}(N/\rho)^{s/3}.
\end{array}
\end{equation}

Firstly, we will prove $F$ and $u$ satisfy the conditions
(\ref{condition2})-(\ref{condition1}) for some $M_j$.  By Lemma
\ref{r0}, we have
\begin{equation}\label{100}
  \| D^ju\|_{ -1/3+n+1, \Omega_{\tilde\rho}} \leq\| D^ju\|_{n+1,\Omega_{\tilde\rho}}
  \leq\frac{C_{7}A^{j-2}}{\tilde\rho^{s(j-3)}}\big((j-3)!\big)^{s} , \quad3\leq j\leq N,
\end{equation}
\begin{equation}\label{+100}
  \|D_v D^ju\|_{ -1/3+n+1, \Omega_{\tilde\rho}}\leq\frac{C_{7}A^{j-2}}
  {\tilde\rho^{s(j-3)}}\big((j-3)!\big)^{s} ,\quad3\leq j\leq N,
\end{equation}
and
\begin{equation}\label{++++100}
\begin{array}{rrr}\| D^ju\|_{ -1/3+n+1, \Omega_{\tilde\rho}}\leq
C_7, \quad 0\leq j\leq2.
\end{array}
\end{equation}
Since $F\in G^s(\bb R^{2n+1}\times \bb R)$, then
\begin{equation}\label{200}
\|(D_{t,x,v}^k\partial_u^lD_p^mF)\big(\cdot,u(\cdot),\nabla_vu(\cdot)\big)\|_{
-1/3+n+1,\Omega}\leq
C_{20}^{k+l}\big((k-3)!\big)^{s}\big((l-3)!\big)^{s}, \quad
k,m+l\geq3.
\end{equation}
Define $M_j, H_0, H_1$ by setting
  \[ H_0=C_7, \quad H_1=A, \quad M_0=C_7 , \indent
     M_j={{\big((j-1)!\big)^{s}}\over{\tilde\rho^{s(j-1)}}}, \quad j\geq1.
  \]
We can choose $A$ large enough such that $H_1=A\geq C_2 H_0$. Then
(\ref{100})-(\ref{200}) can be rewritten
\begin{eqnarray}\label{+a}
  \| D^ju\|_{ -1/3+n+1, \Omega_{\tilde\rho}}\leq H_0,\quad 0\leq j\leq1,
\end{eqnarray}
\begin{eqnarray}\label{+b}
  \| D^ju\|_{ -1/3+n+1, \Omega_{\tilde\rho}}
  \leq H_0H_1^{j-2}M_{j-2},\quad 2\leq j\leq|\alpha|=N,
\end{eqnarray}
\begin{eqnarray}\label{+b'}
  \|D_v D^ju\|_{ -1/3+n+1, \Omega_{\tilde\rho}}
  \leq H_0H_1^{j-2}M_{j-2},\quad 2\leq j\leq|\alpha|=N,
\end{eqnarray}
\begin{eqnarray}\label{+c}
\|(D_{t,x,v}^k\partial_u^lD_p^mF)\|_{ -1/3+n+1,\Omega}\leq
C_{20}^{k+l}M_{k-2}M_{m+l-2},\hspace{0.6cm} k,m+l\geq2.
\end{eqnarray}
For each $j$, note that $s\geq3$ and hence
\begin{equation}\label{328}
\begin{array}{lll}
\frac{j!}{i!(j-i)!}M_iM_{j-i}&=&
\frac{j!}{i(j-i)}\big((i-1)!\big)^{s-1}\big((j-i-1)!\big)^{s-1}
\tilde\rho^{-s(i-1)}
\tilde\rho^{-s(j-i-1)}\\\\
&\leq&(j!)\big((j-2)!\big)^{s-1}\tilde\rho^{-s(j-1)} \\\\
&\leq&\frac{j}{(j-1)^{s-1}}(j-1)!\big((j-1)!\big)^{s-1}\tilde\rho^{-s(j-1)}\\\\
&\leq& M_j.
\end{array}
\end{equation}
Thus $M_j$ satisfy the monotonicity condition (\ref{monotonicity}).
In virtue of (\ref{+a})-(\ref{328}),  using Lemma \ref{+Friedman},
we have
\begin{eqnarray*}
\|\varphi_{\rho,N}D^\alpha[F(\cdot,u(\cdot))]\|_{ -1/3+n+1}&\leq&
C_3H_0H_1^{|\alpha|-2}M_{|\alpha|-2}\\&\leq&
\frac{C_3C_{7}A^{|\alpha|-2}}{\tilde\rho^{s(|\alpha|-3)}}\big((|\alpha|-3)!\big)^{s}\\
&\leq&
\frac{C_{21}A^{|\alpha|-2}}{\rho^{s(|\alpha|-3)}}\big((|\alpha|-3)!\big)^{s}(N/\rho)^{s/3}.
\end{eqnarray*}
This completes the proof of conclusion (\ref {step2}).

\bigbreak {\bf \emph{Step 3.}} Claim
\begin{eqnarray}
\label{step3}\|\mc{L}\varphi_{\rho,N}D^\alpha u\|_{-1/3+n+1}\leq
\frac{C_{23}A^{|\alpha|-2}}{\rho^{s(|\alpha|-3)}}\big((|\alpha|-3)!\big)^{s}(N/\rho)^{s/3}.
\end{eqnarray}

In fact,
\begin{eqnarray*}
\|\mc{L}\varphi_{\rho,N}D^\alpha u\|_{-1/3+n+1} &\leq&
C_{22}\{~\|[\mc{L},~\varphi_{\rho,N}D^\alpha]
u\|_{-1/3+n+1}+\|\varphi_{\rho,N}D^\alpha\mc{L}
u\|_{-1/3+n+1}~\}\\
&=& C_{22}\bigset{\|[\mc{L},~\varphi_{\rho,N}D^\alpha]
u\|_{-1/3+n+1}\\
&&+\norm{\varphi_{\rho,N}D^\alpha[F\big(\cdot,u(\cdot),\nabla_vu(\cdot)\big)]
}_{-1/3+n+1}}.
\end{eqnarray*}
This along with (\ref{step1}), (\ref {step2}) in step 1 and step 2
yields immediately the conclusion (\ref{step3}).

\bigbreak {\bf \emph{Step 4.}} Claim
\begin{equation}\label{step4}
\begin{array}{lll}
\|\varphi_{\rho,N}D^\alpha
u\|_{1/3+n+1}+\|\varphi_{\rho,N}D_vD^\alpha u\|_{1/3-1/3+n+1}\leq
\frac{C_{31}A^{|\alpha|-2}}{\rho^{s(|\alpha|-3)}}\big((|\alpha|-3)!\big)^{s}(N/\rho)^{s/3}.
\end{array}
\end{equation}

In fact, applying the subelliptic estimate (\ref{subestimate}), we
obtain
\begin{eqnarray*}
&&\|\varphi_{\rho,N}D^\alpha u\|_{1/3+n+1} \leq
C_{24}\{~\|\mc{L}\varphi_{\rho,N}D^\alpha
u\|_{-1/3+n+1}+\|\varphi_{\rho,N}D^\alpha u\|_{n+1}~\}.
\end{eqnarray*}
Combining  Lemma \ref{r0} and (\ref{step3}) in Step 3, we have
\begin{eqnarray}\label{N}
\|\varphi_{\rho,N}D^\alpha u\|_{1/3+n+1}\leq
\frac{C_{25}A^{|\alpha|-2}}{\rho^{s(|\alpha|-3)}}\big((|\alpha|-3)!\big)^{s}(N/\rho)^{s/3}.
\end{eqnarray}
Now it remains to treat $\|\varphi_{\rho,N}D_v D^\alpha
u\|_{1/3-1/3+n+1}$, and  $$\|\varphi_{\rho,N}D_v D^\alpha
u\|_{1/3-1/3+n+1}\leq\|D_v \varphi_{\rho,N}D^\alpha
u\|_{n+1}+\|[D_v,~\varphi_{\rho,N}] D^\alpha u\|_{n+1}.$$ Firstly,
we treat the first term on the right. By direct calculation, it
follows that
\begin{eqnarray*}
&&\|D_v\varphi_{\rho,N}D^\alpha
u\|_{n+1}^2\\
&=&{\rm Re}\big(\mc L \varphi_{\rho,N}D^\alpha u,
~a^{-1}\Lambda^{2n+2}\varphi_{\rho,N}D^\alpha u\big)
-{\rm Re}\big(X_0\varphi_{\rho,N}D^\alpha u,~a^{-1}\Lambda^{2n+2}\varphi_{\varepsilon,k\varepsilon}D^\alpha u\big)\\
&=&{\rm Re}\big(\mc L\varphi_{\rho,N}D^\alpha u,
~a^{-1}\Lambda^{2n+2}\varphi_{\rho,N}D^\alpha
u\big)-{1\over2}\big(\varphi_{\rho,N}D^\alpha
u,~[a^{-1}\Lambda^{2n+2},~X_0]\varphi_{\rho,N}D^\alpha
u\big)\\
&&-{1\over2}\big(\varphi_{\rho,N}D^\alpha
u,~[\Lambda^{2n+2},~a^{-1}]X_0\varphi_{\rho,N}D^\alpha
u\big)\\
&\leq&C_{26}\big\{~ \|\mc L\varphi_{\rho,N}D^\alpha
u\|_{-1/3+n+1}^2+\|\varphi_{\rho,N}D^\alpha u\|_{1/3+n+1}^2~\big\}.
\end{eqnarray*}
This along with  (\ref{step3}) and (\ref{N}) shows at once
\begin{eqnarray*}\label{+++}
\|D_v\varphi_{\rho,N}D^\alpha u\|_{r-1/3+n+1}\leq
\frac{C_{27}A^{|\alpha|-2}}{\rho^{s(|\alpha|-3)}}\big((|\alpha|-3)!\big)^{s}(N/\rho)^{s/3}.
\end{eqnarray*}
Moreover Lemma \ref{r0} yields
\begin{eqnarray*}
\|[D_v,~\varphi_{\rho,N}] D^\alpha u\|_{n+1} &\leq& C_{28}(N/\rho)\|
D^\alpha u\|_{n+1, \Omega_{\tilde \rho}}\\&\leq&
\frac{C_{28}C_7A^{|\alpha|-2}}{\tilde\rho^{s(|\alpha|-3)}}\big((|\alpha|-3)!\big)^{s}(N/\rho)^{s/3}\\&\leq&
\frac{C_{29}A^{|\alpha|-2}}{\rho^{s(|\alpha|-3)}}\big((|\alpha|-3)!\big)^{s}(N/\rho)^{s/3}.
\end{eqnarray*}
From the above two inequalities, we have
\begin{eqnarray*}\label{NN}
\|\varphi_{\rho,N}D_vD^\alpha u\|_{1/3+n+1}\leq
\frac{C_{30}A^{|\alpha|-2}}{\rho^{s(|\alpha|-3)}}\big((|\alpha|-3)!\big)^{s}(N/\rho)^{s/3}.
\end{eqnarray*}
This completes the proof of Step 4.

\bigbreak It's clear for any $\rho,~~0<\rho<1,$
\[\|D^\alpha u\|_{1/3+n+1,\Omega_{\rho}}+\|D_vD^\alpha
u\|_{1/3-1/3+n+1,\Omega_{\rho}}\leq\|\varphi_{\rho,N}D^\alpha
u\|_{1/3+n+1}+\|\varphi_{\rho,N}D_vD^\alpha u\|_{1/3-1/3+n+1}.\]
Thus from Step 4, it follows that the conclusion in Lemma \ref{r1/3}
is true for $r=1/3$.

Moreover for any $0< r< 1/3,$ using the interpolation inequality
(\ref{interpolation}), we have
\begin{eqnarray*}
\|D^\alpha u\|_{r+n+1,\Omega_{\rho}}&\leq&
\|\varphi_{\rho,N}D^\alpha u\|_{r+n+1}\\
&\leq&\varepsilon\|\varphi_{\rho,N}D^\alpha
u\|_{1/3+n+1}+\varepsilon^{-r/(1/3-r)}\|\varphi_{\rho,N}D^\alpha
u\|_{n+1}\\
&\leq&\varepsilon\frac{C_{31}A^{|\alpha|-2}}{\rho^{s(|\alpha|-3)}}\big((|\alpha|-3)!\big)^{s}(N/\rho)^{s/3}
+\varepsilon^{-r/(1/3-r)}\frac{C_{32}A^{|\alpha|-2}}{\rho^{s(|\alpha|-3)}}\big((|\alpha|-3)!\big)^{s},
\end{eqnarray*}
Taking $\varepsilon=(N/\rho)^{s(r-1/3)}$, then
\begin{eqnarray*}
\|D^\alpha u\|_{r+n+1,\Omega_{\rho}}
\leq\frac{C_{33}A^{|\alpha|-2}}{\rho^{s(|\alpha|-3)}}\big((|\alpha|-3)!\big)^{s}(N/\rho)^{rs}.
\end{eqnarray*}
Similarly,
\begin{eqnarray*}
\|D_vD^\alpha u\|_{r-1/3+n+1,\Omega_{\rho}}
\leq\frac{C_{34}A^{|\alpha|-2}}{\rho^{s(|\alpha|-3)}}\big((|\alpha|-3)!\big)^{s}(N/\rho)^{rs}.
\end{eqnarray*}
This completes the proof of Lemma \ref{r1/3}.

\bigbreak Inductively,  we have the following
\begin{lem}\label{r1}

For any $r$ with $1/3\leq r\leq2/3$,
\begin{align}
\label{++las}
  \|D^\alpha
  u\|_{r+n+1,\Omega_\rho}+\|D_v D^\alpha
  u\|_{r-1/3+n+1,\Omega_\rho}\leq
  \frac{C_{38}A^{|\alpha|-2}}{\rho^{s(|\alpha|-3)}}\big((|\alpha|-3)!\big)^{s}(N/\rho)^{sr},
  \hspace{0.3cm} \forall~~0<\rho<1.
\end{align}
Moreover, the above inequality still holds for any $r$ with $2/3\leq
r\leq1.$

\end{lem}

\noindent{\bf Proof.} Repeating the proof of Lemma \ref{r1/3}, we
have the truth of (\ref{++las}) for $1/3\leq r\leq 2/3$. The case
$2/3\leq r\leq 1$ is a little different. The conclusion in Step 1 in
the above proof still holds for $r=1$, and  corresponding to Step 2,
we have to make some modification to prove
\begin{equation*}
\begin{array}{lll}
\norm{
\varphi_{\rho,N}D^\alpha[F\big(\cdot,u(\cdot),\nabla_vu(\cdot)\big)]}_{
1/3+n+1}\leq
\frac{C_{36}A^{|\alpha|-2}}{\rho^{s(|\alpha|-3)}}\big((|\alpha|-3)!\big)^{s}(N/\rho)^{s}.
\end{array}
\end{equation*}
From the truth of (\ref{++las})  for $1/3\leq r\leq 2/3$, it follows
\[\| D^ju\|_{ 1/3+n+1, \Omega_{\tilde\rho}}
  \leq\frac{C_{37}A^{j-2}}{\tilde\rho^{s(j-3)}}\big((j-3)!\big)^{s}(j/\tilde\rho)^{s/3},
  \quad 3\leq j\leq N,
\]
\[\|D_v D^ju\|_{ 1/3+n+1, \Omega_{\tilde\rho}}
  \leq\|D_v D^ju\|_{2/3- 1/3+n+1, \Omega_{\tilde\rho}}
  \leq\frac{C_{37}A^{j-2}}{\tilde\rho^{s(j-3)}}\big((j-3)!\big)^{s}
  (j/\tilde\rho)^{2s/3},\quad 3\leq j\leq N,
\]
and
\begin{equation*}
\begin{array}{rrr}\| D^ju\|_{ 1/3+n+1, \Omega_{\tilde\rho}}\leq
C_{37}, \hspace{0.6cm}0\leq j\leq2,
\end{array}
\end{equation*}
Hence we need define a new sequence  $\bar M_j$ by setting
\[\bar M_0=C_{37},\indent \bar M_j={{\big((j-1)!\big)^{s}}
  \over{\tilde\rho^{s(j-1)}}}\big((j+2)/\tilde\rho\big)^{2s/3},\quad j\geq1.\]
  For each $j$, note that $s\geq3$ and hence
direct computation deduces that for $0<i<j,$
\begin{equation*}
\begin{array}{lll}
  \frac{j!}{i!(j-i)!}\bar M_i\bar M_{j-i}
  &=&\frac{j!}{i(j-i)}\big((i-1)!\big)^{s-1}\big((j-i-1)!\big)^{s-1}\\\\
    &&\times (i+2)^{2s/3}(j-i+2)^{2s/3}\tilde\rho^{-s(j-2)}\tilde\rho^{-4s/3}\\\\
  &\leq&4(j!)\big((j-2)!\big)^{s-1}(j+2)^{2s/3-1}(j+1)^{2s/3-1}
    \tilde\rho^{-s(j-1)}\tilde\rho^{-2s/3}\tilde\rho^{s-2s/3}\\\\
  &\leq&\frac{4j(j+1)^{2s/3-1}}{(j-1)^{s-1}}(j-1)!\big((j-1)!\big)^{s-1}
    \tilde\rho^{-s(j-1)}\big((j+2)/\tilde\rho\big)^{2s/3}\\\\
  &\leq& C_{39}\bar M_j.
\end{array}
\end{equation*}
In the last inequality we used the fact $s-1\geq2s/3$. Thus $\bar
M_j$ satisfy the monotonicity condition (\ref{monotonicity}). Now
the left is entirely  similar to the proof of Lemma \ref{r1/3}. And
thus (\ref{++las}) holds for $r=1$ and hence for $2/3\leq r\leq 1$
by interpolation inequality (\ref{interpolation}). This completes
the proof of Lemma \ref{r1}.

\bigbreak Recall $C_7,~C_{35}$ and $C_{35}$ are the constants
appearing in Lemma \ref{r0}, Lemma \ref{r1/3} and Lemma \ref{r1}.
Now taking $A$ large enough such that $A\geq\max\{C_7,
C_{35},C_{38}\}$, and then by the above three Lemmas we have the
truth of $(E)_{r,N}$ for any $r\in[0,~1].$ This complete the proof
of Proposition \ref{prp4'}.


\begin{thebibliography}{99}

\bibitem{aumxy} R. Alexandre, S. Ukai, Y. Morimoto, C.-J. Xu, T. Yang,
Uncertainty principle and regularity for Boltzmann equation,


\bibitem{al-1} R. Alexandre, L. Desvillettes, C. Villani, B. Wennberg,
Entropy  dissipation and long-range interactions, {\em Arch.
Rational Mech. Anal.} {\bf 152} (2000) 327-355.

\bibitem{al-2} R. Alexandre, M. Safadi, Littlewood Paley decomposition and
regularity issues in Boltzmann equation homogeneous equations. I.
Non cutoff and Maxwell cases, M3AM (2005) 8-15.


\bibitem{Bou}\label{Bou}F. Bouchut,~
Hypoelliptic regularity in kinetic equations. {\sl J. Math. Pure
Appl.} {\bf 81} (2002), 1135-1159.


\bibitem{CR}\label{CR} Chen Hua,  L.Rodino,
General theory of PDE and Gevrey class. {\sl General theory of
partial differential equations and microlocal analysis}(Trieste
1995), Pitman Res. Notes in Math. Ser., {\bf 349}, Longman, Harlow,
6-81, (1996).

\bibitem{CHR} Chen Hua,  L.Rodino,
Paradifferential calculus in Gevrey class . {\sl J. Math. Kyoto
Univ.} {\bf 41}, (2001), 1-31.


\bibitem{DZ}\label{DZ}M.Derridj,   C.Zuily,
Sur la r\'{e}gularit\'{e} Gevrey des op\'{e}rateurs de
H\"{o}rmander. {\sl J.Math.Pures et Appl}. {\bf 52} (1973), 309-336.

\bibitem{desv-wen1} L. Desvillettes, B. Wennberg, Smoothness of the solution
of the spatially homogeneous Boltzmann equation without cutoff. {\em
Comm. Partial Differential Equations} {\bf 29} (2004), no. 1-2,
133--155.

\bibitem{Dur} \label{Dur}M.Durand,
R\'{e}gularit\'{e} Gevrey d'une classe d'op\'{e}rateurs
hypo-elliptiques. {\sl J.Math.Pures et Appl}. {\bf 57} (1978),
323-360.

\bibitem{Friedman} A.Friedman,
On the Regularity of the solutions of Non-linear Elliptic and
Parabolic Systems of Partial Differential Equations. {\sl J. Math.
Mech}. {\bf 7} (1958), 43-59.


\bibitem{helffer-nier} B. Helffer, F. Nier, Hypoelliptic estimates and
spectral theory for Fokker-Planck operators and Witten Laplacians.
{\em Lecture Notes in Mathematics}, {\bf 1862} Springer-Verlag,
Berlin, 2005.


\bibitem{herau-nier}F. H\'erau, F. Nier, Isotropic hypoellipticity and
trend to equilibrium for the Fokker-Planck equation with a
high-degree potential. {\em Arch. Ration. Mech. Anal.} {\bf 171}
(2004), no. 2, 151--218.

\bibitem{Hor}\label{Hor}L.H\"{o}rmander,
Hypoelliptic second order differential equations. {\sl Acta Math.}
{\bf 119} (1967), 147-171.


\bibitem{Kohn}  \label{Kohn}J.Kohn,
Lectures on degenerate elliptic problems. Psedodifferential
operators with applications, C.I.M.E., Bressanone 1977,
89-151(1978).


\bibitem{MX}\label{MX}Y.Morimoto, C.-J. Xu,
Hypoellipticity for a class of kinetic equations. to appear at J.
Math. Kyoto U.

\bibitem{mo-xu}Y. Morimoto and C.-J. Xu, Logarithmic Sobolev inequality
and semi-linear Dirichlet problems for infinitely degenerate
elliptic operators, {\it Ast\'erisque} {\bf 284} (2003), 245--264.

\bibitem{MUXY1} Y. Morimoto, S. Ukai, C.-J. Xu, T. Yang,  Regularity
of solutions to the spatially homogeneous Boltzmann equation without
Angular cutoff, preprint.


\bibitem{Ro} \label{Ro} L.Rodino,
{\sl Linear partial differential operators in Gevrey class.}  World
Scientific, Singapore, 1993.


\bibitem {R-S}\label{RE}L.P.Rothschild,  E.M.Stein,
Hypoelliptic differential operators and nilpotent groups.  {\sl
Acta.Math}. {\bf 137} (1977), 248-315.


\bibitem {Tr}\label{Treves} F.Treves,
{\sl Introduction to Pseudodifferential and Fourier Integral
Operators.} Plenum, New York, 1980.


\bibitem{ukai}S. Ukai, Local solutions in Gevrey classes to the nonlinear Boltzmann equation
without cutoff, {\em Japan J. Appl. Math.}{\bf 1}(1984), no. 1,
141--156.


\bibitem{villani}C. Villani, On a new class of weak solutions to the spatially
homogeneous Boltzmann and Landau equations, {\em Arc. Rational Mech.
Anal., } {\bf 143}, 273--307, (1998).


\bibitem{Xu} C.-J. Xu,  Nonlinear microlocal analysis. {\sl General theory
of partial differential equations and microlocal analysis}(Trieste
1995), Pitman Res. Notes in Math. Ser., {\bf 349}, Longman, Harlow,
155-182, (1996).

\end{thebibliography}
\end{document}